\documentclass{article}
\usepackage{amsmath,ifart,ifmath,palatino,mathpazo,graphics}
\DeclareMathOperator{\kbes}{K} 
\newcommand{\dom}{\Omega} 
\newcommand{\eNrg}{\mathcal E} 
\begin{document}

\iftitle{A study of blow-ups in the Keller-Segel model of chemotaxis}

\ifauthor{Ibrahim Fatkullin}{Department of Mathematics, University of Arizona}

\ifabstract{We study the Keller-Segel model of chemotaxis and develop a composite particle-grid numerical method with adaptive time stepping which allows us to resolve and propagate singular solutions. We compare the numerical findings (in two dimensions) with analytical predictions regarding formation and interaction of singularities obtained through analysis of the stochastic differential equations associated with the model.}


\vspace{2ex}
{\hfill\small\today}

\vspace{-2ex}
\section{Introduction}

In this work we discuss a numerical method for studying singularities (blow-ups) in the Keller-Segel model (K-S) of chemotaxis \cite{Kell70} and in the systems with similar mathematical structure, such as the McKean-Vlasov (McK-V) model of charged interacting particles \cite{Vill03}. (Patlack derived similar equations nearly two decades prior to Keller and Segel \cite{Pat53}.) We present illustrative numerical examples and provide analytical insights into formation and interaction of singularities from the point of view of the associated stochastic processes. We concentrate on the details of numerical modeling and several specific aspects of mathematical analysis. An interested reader should consult the exhaustive reviews for discussions of physical and biochemical phenomena involved in chemotaxis \cite{Hor03,Hor04,Tin08i,Tin08ii}.

Numerical treatment of equations with singular solutions is challenging because it requires a reliable and efficient approximation of singular functions or even distributions. For example, the conventional finite-difference or spectral methods are not fit for dealing with $\delta$-functions arising in the K-S model. While most finite-element and discontinuous Galerkin methods are efficient and accurate as long as the solutions remain sufficiently smooth, in the vicinity of singularities they struggle with positivity issues, unphysical oscillations, and related problems \cite{Fil06,Epsh08,Epsh09,Epsh09a,Stre10a,Stre10b}. (A conservative upwind finite-element method introduced by Saito for the elliptic K-S model overcomes some of these issues \cite{Sai07}.) Generally, unless some regularizations are introduced, these methods cannot propagate solutions past the moment of blow-up altogether.
The method presented in this work avoids such difficulties because the singular field is modeled by an ensemble of interacting particles. The singularities manifest themselves as particle aggregates, harmless from the numerical standpoint. Essentially, this approach to chemotaxis reverses the continuous PDE description of ensembles of interacting particles.
Ha\v skovec and Schmeiser have recently introduced a method utilizing similar ideas \cite{Has10,Has11}. Their method, however, is exclusively particle-based and is limited to one specific variety of the elliptic K-S model in the entire plane. Our method employs the grid and particle representations of the fields simultaneously and is suitable for both elliptic and parabolic K-S equations in arbitrary domains. It is also applicable to a wider range of kinetic problems. (Other technical differences are discussed where appropriate.) Conceptually, our method is more in the spirit of the Particle In Cell (PIC) methods widely used, e.g., in plasma physics \cite{Iva80,Cot86,Bird91,Vic91,Allen94,Hav96,Wall00}.
Another interesting numerical method has been recently implemented for the one-dimensional McKean-Vlasov system \cite{Bla08}. It uses that the McK-V equation describes steepest descent dynamics in the Wasserstein-2 space of probability measures and is suitable for studying singular (measure-valued) solutions. 

The usefulness of the particle approach extends beyond numerical methods. Analysis of the stochastic processes underlying the K-S PDEs yields insights into formation and interaction of singularities. In particular, we derive an expression for the critical mass required to create or sustain a singularity and relate the non-uniqueness of solutions of the K-S equations to the underlying diffusion process. Even though our treatment is rather informal, it sheds new light on the mechanism of blow-ups.
Some of the results presented in this paper have been proven using different methods and under various assumptions. Much effort has been directed towards obtaining the exact value of the critical mass \cite{Dol03,Bla06,Calv06,Chav07,Bil09}. Vel\'azquez introduced measure-valued solutions and extended the K-S equations beyond the point of blow-up \cite{Vel04a,Vel04b}, see also \cite{Dol09}. These works together with the study of the K-S equations as a hydrodynamic limit of interacting particle systems by Stevens \cite{Stev00} link the particle and PDE descriptions of chemotaxis.
We also mention a rigorous analysis of phase transitions in the McK-V system \cite{Chay09} (it only covers sufficiently regular interaction kernels not relevant to the K-S model) and an investigation of singularity formation in several related aggregation models \cite{Bert09,Bert10,Carri10} (these models do not have diffusive terms and correspond to one special case of the K-S model).

\paragraph{The Keller-Segel model} is prescribed by the following system of PDEs:
\begin{subequations}\label{eq:KS}
\begin{align}
	\pd_t\rho(\vv x,t)\,&=\,\grd\cdot(\mu\grd\rho\,-\,\chi\,\rho\grd c)\label{eq:KSrho},\\
  	\alpha\,\pd_tc(\vv x,t)\,&=\,\Delta c\,-\,k^2c\,+\,\rho.\label{eq:KSc}
\end{align}
\end{subequations}
The function $\rho(\vv x)$ is the density of active particles (bacteria), $c(\vv x)$ is the concentration of chemoattractant. For numerical simulations we use Neumann (no flux) boundary conditions, 
\begin{equation}
	\pd_{\vv n}c(\vv x)=\pd_{\vv n}\rho(\vv x)=0\quad\text{for}\quad\vv x\in\pd\dom
\end{equation} 
in a two-dimensional square domain $\dom=(0,L)^2$.
We assume that $\rho(\vv x)$ integrates to $M$ over the entire domain; $M$ is the total mass of the bacteria. By rescaling the equations, we can always achieve that $\alpha=0$ or $\alpha=1$. In the former case the model is called {\em elliptic,} in the latter --- {\em parabolic.} The parameters $\mu$, $\chi$, and $k$ are constants. (The numerical method, however, may be readily extended to a more general class of equations, see Section~\ref{sec:num} for details.)

Physically, these equations describe the following phenomenon. The chemoattractant spreads diffusively and decays with rate $k^2$; it is also produced by the bacteria with rate $1$. (In the elliptic case these rates are infinite, i.e., the chemoattractant ``thermalizes'' infinitely fast.) The bacteria diffuse with mobility $\mu$ and also drift in the direction of the gradient of concentration of the chemoattractant with velocity $\chi|\grd c|$; $\chi$ is called chemosensitivity.

In this work we are interested in singular solutions to the K-S equations. A typical situation is illustrated in Figure~\ref{fig:blowups} where a snapshot of the concentration field $c(\vv x)$ is displayed. The peaks correspond to $\delta$-function-type singularities of the particle density $\rho(\vv x)$. (Displaying bacteria density itself is not very illustrative due to its highly singular nature.) As mentioned above, the principal challenge for numerical simulations is in approximating the fields $\rho(\vv x)$ and $c(\vv x)$ as they become unbounded. Our idea is to use a particle method for the evolution of $\rho(\vv x)$, so that its singularities manifest themselves as harmless particle aggregates. For the propagation of $c(\vv x)$, we use a second order implicit finite-difference scheme because of its simplicity and excellent stability properties. Even though $c(\vv x)$ also becomes unbounded analytically (it remains large but bounded for a given discretization), its weak logarithmic singularities do not cause any problems for sufficiently stable schemes. In general, one may use a different solver for propagation of the concentration field (including a particle-based solver as well).

\begin{figure}[h]
\scalebox{0.4}{\includegraphics{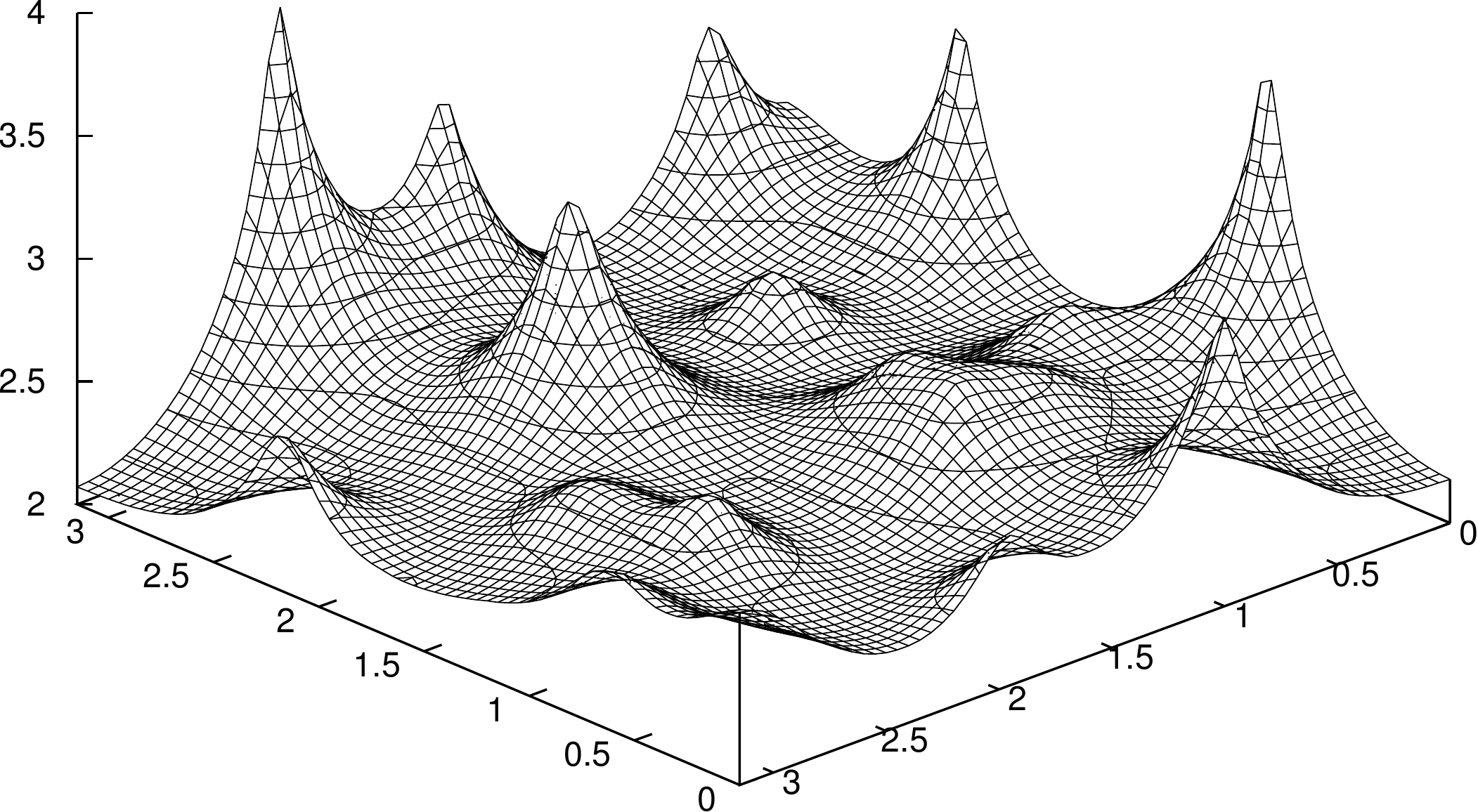}}\hfill
\parbox[b]{2.25in}{\caption{\label{fig:blowups}\small A snapshot ($t=10$) of the concentration field $c(\vv x)$ illustrating multiple blow-ups in the parabolic model. The parameters are as follows: $\alpha=k=1$, $\chi=0.1$, $\mu=0.005$, $M=25$. The spatial domain, $3.2\times3.2$, is discretized using the mesh size $\Delta x=0.05$; the time step for propagation of the concentration field is $\Delta t=0.1$; 4096 particles are used to approximate $\rho(\vv x)$ (see Section~\ref{sec:num} for details of the numerical scheme). Some random initial data is used at $t=0$.}}
\end{figure}

\paragraph{The particle-based description} utilizes that the PDE (\ref{eq:KSrho}) is a Fokker-Planck equation for the evolution of the probability density of the bacteria. Suppose first, that we are given some time-dependent concentration field $c(\vv x,t)$ and want to solve equation (\ref{eq:KSrho}). Let $\vv X^{(n)}_t$ be i.~i.~d. random variables distributed according to the initial condition $\rho(\vv x,0)$ at time $t=0$ and satisfying the following stochastic differential equations (with reflecting boundary at $\pd\dom$) for $t>0$:
\begin{equation}\label{eq:KS_SDE}
	\md \vv X^{(n)}_t\,=\,\chi\grd c\big(\vv X^{(n)}_t,t\big)\md t\,+\,\sqrt{2\mu}\,\md\vv W^{(n)}_t.
\end{equation}
Then, using the strong law of large numbers and that (\ref{eq:KSrho}) is a Fokker-Planck equation for the SDEs (\ref{eq:KS_SDE}), we can see that as $N\to\infty$, the empirical probability densities,
\begin{equation}\label{eq:ePDF}
	P_N(\vv x,t)\,=\,\frac{M}{N}\sum_{n=1}^N\delta\big(\vv x-\vv X^{(n)}_t\big),
\end{equation}
converge almost surely to $\rho(\vv x,t)$. (This convergence is in the sense of measures, point-wise in $t$; a stronger convergence can be established with more effort.)
Thus we can simulate equations (\ref{eq:KS_SDE}) for the particles, and approximate the density $\rho(\vv x,t)$ using $P_N(\vv x,t)$ whenever required.

Rigorously justifying this approach to the full K-S system (\ref{eq:KS}) is more challenging because in the full system the concentration field $c(\vv x,t)$ is a functional of $\rho(\vv x, s)$, $s\leq t$, rather than some a priori prescribed field. This implies that equations (\ref{eq:KS_SDE}) become coupled and we end up with a system of stochastic interacting particles (with memory in the parabolic case). Some progress in this direction, nevertheless, has been made.  Stevens derived the K-S equations as limit dynamics for interacting stochastic particle systems via smoothing and rescaling the interaction potentials in a particular fashion \cite{Stev00}. Ha\v skovec and Schmeiser proved convergence of their particle method under reasonable assumptions via analysis of the associated BBGKY hierarchy \cite{Has11}.


\section{Numerical method}\label{sec:num}
As mentioned in the introduction, the principal idea behind our numerical scheme is to employ a particle method for the evolution of the particle density field $\rho(\vv x,t)$. The current implementation utilizes an implicit second order finite-difference scheme for simulating equation (\ref{eq:KSc}) and an explicit Euler-Maruyama scheme for the stochastic particle dynamics (\ref{eq:KS_SDE}). We discretize the computational (rectangular) domain with grid size $\Delta x$ and propagate the concentration field $C_{ij}$ using the time step $\Delta t$.  The particle density field $P_{ij}$ is reconstructed from the particle locations. We evolve the particles using adaptive time steps which may be smaller than $\Delta t$; this is needed for stability reasons (see below).

\paragraph{The particle dynamics}
is simulated using the forward Euler-Maruyama scheme,
\begin{equation}
  \vv X^{(n)}(t+\Delta\tau)\,=\,\vv X^{(n)}(t)\,+\,\chi\grd c\big(\vv X^{(n)}(t)\big)\Delta\tau\,+\,\sqrt{2\mu\,\Delta\tau}\,\vv N(0,1),
\end{equation}
where $N(0,1)$ is a standard Gaussian random variable with mean 0 and variance 1. The gradient field $\grd c(\vv x)$ is approximated in two steps. First, we construct the gradient fields $CX$ and $CY$ using the second order approximation,
\begin{equation}
  CX_{ij}\,=\,\frac{1}{2\Delta x}\Big[C_{i+1,j}(t)\,-\,C_{i-1,j}(t)\Big],\qquad
  CY_{ij}\,=\,\frac{1}{2\Delta x}\Big[C_{i,j+1}(t)\,-\,C_{i,j-1}(t)\Big].
\end{equation}
Then we approximate $\grd c(\vv x)$ via bilinear interpolation using the values of $CX$ and $CY$ at the four nearest grid points. Because $|\grd c(\vv x)|$ becomes unbounded (very large numerically) in the vicinity of blow-ups, it is essential to choose the time step $\Delta\tau$ adaptively, see Figure~\ref{fig:adapt} and its description for details. Consequently, if needed, we subdivide the timestep $\Delta t$ into smaller intervals $\Delta\tau^{(n)}$ so that $|\grd c\big(\vv X^{(n)}\big)|\Delta\tau^{(n)}<\Delta x$. (Each particle is simulated independently of others, i.e, $\Delta\tau^{(n)}$ are chosen for each particle individually.) A failure to satisfy this condition leads to overshooting in the vicinity of blow-ups, which in turn causes artificial oscillations and damping. Finally, we enforce the no-flux Neumann  boundary conditions for $\rho(\vv x)$ by reflecting the particles escaping from the spatial domain back into it.

\begin{figure}
\scalebox{0.7}{\includegraphics{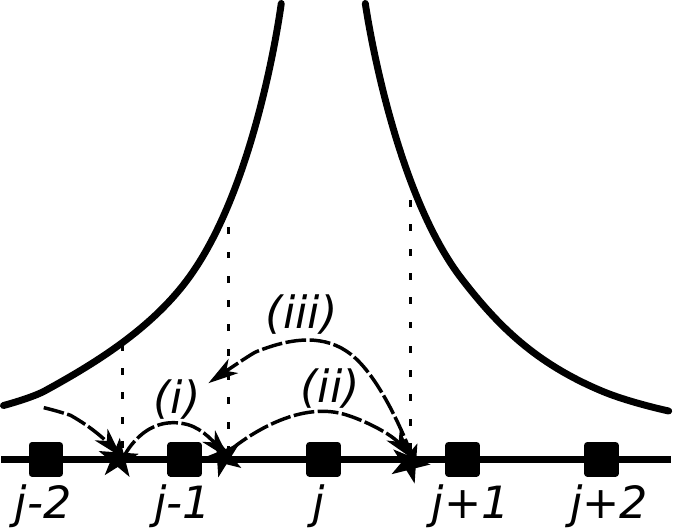}}\hfill
\parbox[b]{3.75in}{\caption{\label{fig:adapt}\small Adaptive time stepping for particle dynamics. In this example a blow-up of the concentration field occurs near a grid point with index $j$. While some particular time step is appropriate for propagating the particles away from the blow-up (i), it may lead to overshooting in its vicinity (ii, iii). Thus the particles incorrectly accumulate near the grid points $j\pm1$ rather than $j$. One can avoid this phenomenon reducing the time step for the particle evolution whenever the gradient of the concentration field becomes too large. Practically, it is sufficient to adjust the time step so that the expected length of the particle jump does not exceed the mesh size $\Delta x$.}}
\end{figure}

\paragraph{The particle density field $P_{ij}(t)$} is reconstructed from the particle locations as explained in the description to Figure~\ref{fig:scheme}: each particle contributes fractions of its weight to the four nearest grid points according to the bilinear interpolation rules.  This approximation preserves the first moment of the corresponding distribution. It is important to use this kind of moment-preserving approximations, e.g., a simple bin-counting (assigning all weight of a given particle to the nearest grid point) is unsatisfactory. The reason is that this creates an artificial flux towards the grid points: the particles feel their own potentials which become artificially aligned to the nearest grid points. This flux is sufficient to pin the singularities to grid points and disrupt such phenomena as logarithmically-weak interaction of singularities with each other and with the boundary of the domain.

\paragraph{The concentration field} is propagated according to the implicit second order finite-difference\linebreak scheme,
\begin{equation}
  \frac{\alpha}{\Delta t}\,\Big[C_{ij}(t+\Delta t)\,-\,C_{ij}(t)\Big]\,=\,\frac{1}{\Delta x^2}\,D^{(2)}_{ij}C(t+\Delta t)\,-\,k^2C_{ij}(t+\Delta t)\,+\,P_{ij}(t),
\end{equation}
where $D^{(2)}$ is the standard second difference operator. It is essential to use an implicit scheme for stability reasons because the concentration field acquires logarithmic singularities. This scheme is quite efficient: even though it requires solving a system of linear equations, the matrix of the finite difference operator is symmetric and banded, and may be Cholesky-factorized  before the actual computations. The extra computational cost of solving this linear system is then alleviated by relaxation of the $\Delta t\sim\Delta x^2$ stability constraint of explicit schemes.

\begin{figure}
\scalebox{0.7}{\includegraphics{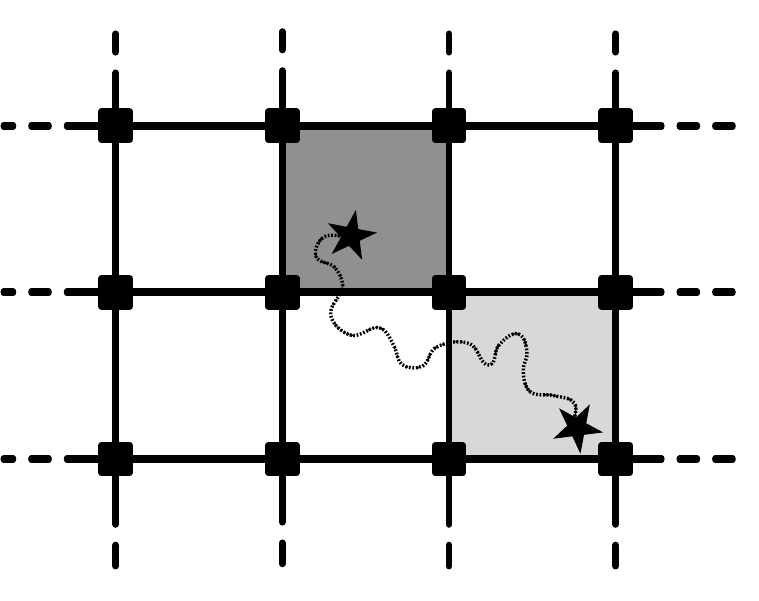}}\hfill
\parbox[b]{3.5in}{\caption{\label{fig:scheme}\small Schematic representation of the numerical scheme. The concentration field $c(\vv x)$ is sampled on a uniform grid ($\blacksquare$), while the particles ($\bigstar$) move without any restraints. The particle density $\rho(\vv x)$ is reconstructed via bilinear interpolation: each particle contributes fractions of its mass to the four nearest grid points proportionally to the relative distances from the latter, e.g., a particle located at $(x,y)$ within a unit-square grid cell contributes weights proportional to $xy$, $x(1-y)$, $y(1-x)$, and $(1-x)(1-y)$ to the respective grid points.}\rule{0pt}{2ex}}
\end{figure}

\newpage
The presented numerical method may be generalized to other similar systems of the form
\begin{subequations}\label{eq:KSgen}
\begin{align}
  \pd_t\rho(\vv x,t)\,&=\,\grd\cdot\big(\mu(\vv x,\rho,c)\grd\rho\,-\,\rho\,\vv v(\vv x,\rho,c)\big)\label{eq:KSrhogen},\\
  \alpha\,\pd_tc(\vv x,t)\,&=\,\grd\cdot\big(G(\vv x)\grd c\big)\,+\,F(\vv x,\rho,c),\label{eq:KScgen}
\end{align}
\end{subequations}
and to more complicated domains. This makes it a useful tool for studying various aggregation phenomena for which the conventional finite-difference or finite-elements schemes are inefficient. From the numerical analysis point of view, we do not currently have any estimates of the convergence rates for this scheme. Partially, this problem is complicated due to the singular nature of solutions to the K-S equations. In particular, once the blow-ups are formed, one must interpret the K-S PDEs in some proper sense, and while our numerical scheme is able to propagate the solutions for all times (i.e., it regularizes the PDEs in some fashion), relation of such numerical solutions to a particular analytical regularization must be carefully investigated. 
On the other hand, we could speculate that our method exhibits the typical for this kind of schemes errors of the order $\bO(\sqrt{\Delta t})$, $\bO(1/\sqrt{N})$, and $\bO(\Delta x^2)$ for regular solutions. Such analysis, however, is not the purpose of this paper and we verify our numerical results by ``standard means,'' i.e., compare the results with simulations employing half the mesh size, double the number of particles, etc.


\section{Formation and interaction of singularities}
We now concentrate on the elliptic ($\alpha=0$) case, though some of the reasoning is equally applicable to the parabolic case as well. Generally, one can distinguish the {\em soft\,} and {\em hard\,} blow-ups: in the first case the fields $\rho(\vv x)$ and $c(\vv x)$ become unbounded but $\rho(\vv x)$ does not acquire atomic components, while in the second case the particle density field becomes a bona fide distribution. Here we concentrate on the hard blow-ups which have a more apparent physical meaning: atomic components of the particle density correspond to accumulation  of particles at some point locations on the given physical length-scale.

It is convenient to study the elliptic model in a slightly different form, as a McKean-Vlasov system. First, we solve equation (\ref{eq:KSc}) for $c(\vv x)$:
\begin{equation}\label{eq:cconv}
  c(\vv x)\,=\,-\int_\dom V_d(\vv x,\vv y)\,\rho(\vv y)\md\vv y\,=\,-\,V_d*\rho(\vv x).
\end{equation}
Here $V_d(\vv x,\vv y)$ is the Green's function for $\Delta-k^2$ in $\dom$, which is assumed to be $d$-dimensional. The entire-space Green's functions (the fundamental solutions)  whose singular part is identical to those in bounded domains are known explicitly:
\begin{subequations}\label{eq:FS}
\begin{align}
	V_1(\vv x,\vv y)\,=&\,-\frac{\me^{-k|\vv x-\vv y|}}{2k}&(d=1);\\
	V_2(\vv x,\vv y)\,=&\,-\frac{\kbes_0(k|\vv x-\vv y|)}{2\pi}&(d=2);\label{eq:fs_2d}\\
	V_3(\vv x,\vv y)\,=&\,-\frac{\me^{-k|\vv x-\vv y|}}{4\pi|\vv x-\vv y|}&(d=3).
\end{align}
\end{subequations}
Here $\kbes_0(\cdot)$ is the modified Bessel function of second kind. Substituting expression (\ref{eq:cconv}) for the concentration field $c(\vv x)$ into equation (\ref{eq:KSrho}) we get a closed integro-differential equation for the particle density field,
\begin{equation}\label{eq:KS_ell_rho}
  \pd_t\rho(\vv x,t)\,=\,\grd\cdot(\mu\grd\rho\,+\,\chi\,\rho\grd V_d*\rho)\,=\,\grd\cdot\left[\rho\grd\VD{\eNrg}{\rho}\right].
\end{equation}
Here the energy functional is given by
\begin{equation}\label{eq:nrg}
  \eNrg(\rho)\,=\,\mu\int_\dom\rho(\vv x)\ln\rho(\vv x)\md\vv x\,+\,\frac{\chi}{2}\iint_\dom\rho(\vv x)V_d(\vv x,\vv y)\rho(\vv y)\md\vv x\md\vv y.
\end{equation}
One of the features of the dynamics (\ref{eq:KS_ell_rho}) is that the functional $\eNrg(\rho)$ is non-increasing:
\begin{equation}
  \PD{\eNrg}{t}\,=\,-\int_{\dom}\left|\grd\VD{\eNrg}{\rho}\right|^2\rho(\vv x,t)\md\vv x\,\leq\,0.
\end{equation}
When $\rho(x)$ becomes singular with respect to Lebesgue measure, the first ({\em entropic}) term in (\ref{eq:nrg}) tends to infinity, which must be compensated by the second ({\em interaction}) term. If $d=1$, the Green's function $V_1(\vv x)$ is bounded below and so is the interaction term, thus the blow-ups are not possible. This is no longer the case if $d\geq 2$ and so the blow-ups are permissible. A more delicate analysis is needed to understand how exactly they are formed.

From standpoint of the stochastic process underlying the Fokker-Planck PDE (\ref{eq:KSrho}), the atomic components of $\rho(\vv x)$ correspond to particles aggregating at point locations. Thus the traps  which do not allow the particles to escape must be created. Let us first investigate when a particle aggregate creates a trap for a single particle diffusing in its field. Suppose $\rho(\vv x)=M\delta(\vv x)$, the corresponding concentration field is then proportional to the fundamental solution (\ref{eq:FS}).  Equation for the radial component of the location of a single particle diffusing in the field created by $\rho$ is equivalent in law to the following SDE:
\begin{equation}
  \md r_t\,=\,\Big[\chi c^\prime\big(r_t\big)\,+\,(d-1)\frac{\mu}{r_t}\Big]\md t\,+\,\sqrt{2\mu}\,\md W_t.
\end{equation}
For $d=1$ we get
\begin{equation}
  \md r_t\,=\,-\frac{\chi M}{2}\exp(-kr_t)\md t\,+\,\sqrt{2\mu}\,\md W_t.
\end{equation}
This equation is sufficiently regular and does not allow for existence of a trap for any value of $M$, so consistently with energy considerations, no blow-ups are possible in this case. For $d\geq 3$, a similar reasoning shows that a trap of an arbitrarily small mass may exist. So the interesting case is $d=2$.

For $d=2$, using that in the leading order, as $r\to0$, $\kbes_0(r)\sim-\ln(r)$, we obtain the following equation:
\begin{equation}\label{eq:bessel}
  \md r_t\,=\,\Big(\mu\,-\,\frac{\chi M}{2\pi}\Big)\frac{\md t}{r_t}\,+\,\sqrt{2\mu}\,\md W_t.
\end{equation}
This is a well-known Bessel process \cite{Kar81}, it's behavior near the boundary at $r=0$ depends strongly on the value of $M$. In particular, when
\begin{equation}
  M\,\geq\,\frac{4\pi\mu}{\chi}\,=\,M_c^*,
\end{equation}
the origin is the so-called {\em exit boundary,} a trap, i.e., the particles reach it in finite time, and may never escape back to the domain $r>0$. For our problem the implication is that the total mass of at least $M_c^*$ is needed for existence of a stable singularity in the particle density field. For $M\in(0,M^*_c)$, the SDE (\ref{eq:bessel}) has the so-called {\em regular boundary} at $r=0$. In this case the particles may reach the origin, but are not necessarily trapped there, i.e., may also leave it according to some rules which must be prescribed in addition to the SDE itself. This implies that the original K-S PDEs alone are not sufficient to describe the blow-up dynamics and the exchange of mass between the regular and singular components of the particle density when the mass of the singularity is less than $M_c^*$.

The critical mass $M_c^*$ is smallest mass such that a singularity does not shed its mass but only absorbs particles from the smooth component. It is, however, twice less than the mass required to create a singularity from smooth initial data. Indeed, consider the $N$-particle SDEs approximating the elliptic K-S equations:
\begin{equation}\label{eq:KS_SDE_n}
  \md \vv X^{(n)}_t\,=\,-\frac{\chi M}{N}\Pd{}{\vv X^{(n)}}\sum_{i=1,\;i\neq n}^N
  V_d\big(\vv X^{(n)}_t\!\!,\,\vv X^{(i)}_t\big)\md t\,+\,\sqrt{2\mu}\,\md\vv W^{(n)}_t.
\end{equation}
%
%
%
%
%
%
%
%
%
To study the aggregation dynamics of the particle system, consider the variance of the empirical probability density of the bacteria (\ref{eq:ePDF}),
\begin{equation}\label{eq:R}
	R_t^2\,=\,\frac{1}{\,2N^2}\sum_{i,\,j=1}^{N}\big|\vv X^{(i)}_t-\,\vv X^{(j)}_t\big|^2;
\end{equation}
$R_t$ will be referred to as the\, {\em radius\,} of the system. We concentrate on the problem in the entire space with $k=0$; this corresponds to $V_2(\vv x,\vv y)=1/(2\pi)\ln|\vv x-\vv y|$. Analysis for the other values of $k$, or in the bounded domains is similar (although more technical) because the singularities of $V_2(\vv x,\vv y)$  are the same in the leading order. In this case, however, there exists an explicit closed equation governing the evolution of $R_t$:
\begin{equation}\label{eq:SDE_N}
	\md R_t\;=\,\Bigg[2\mu\,\bigg(1-\frac{3}{2N}\bigg)-\,\frac{\chi M}{4\pi}\bigg(1-\frac{1}{N}\bigg)\Bigg]\;\frac{\md t}{R_t}\;+\;\sqrt{\frac{2\mu}{N}\;}\md W_t.
\end{equation}
The driving stochastic process $W_t$ obeys the following SDE:
\begin{equation}
	\md W_t\,=\,\frac{1}{N^{3/2}\,R_t}\;\sum_{i,\,j=1}^{N}\Big(\vv X^{(i)}_t-\,\vv X^{(j)}_t\Big)\cdot\,\md\vv W^{(i)}_t,
\end{equation}
and can be shown to be a Wiener process using the L\'evy characterization of Brownian motion. We can immediately see that in the limit as $N\to\infty$, if
\begin{equation}
 	M\,\neq\,\frac{8\pi\mu}{\chi}\,=\,M_c,
\end{equation}
the evolution of the radius obeys a well-known deterministic equation, 
\begin{equation}
	\dot R_t\,=\,-\frac{\gamma}{2R_t},\qquad
	\gamma\,=\,4\mu(M/M_c-1).
\end{equation}
Hence $R^2_t=R^2_0-\gamma t\,$, and for supercritical masses ($M>M_c$) the blow-up occurs at time $T=R_0^2/\gamma$. It is independent of the initial distribution of mass except through $R_0$. If $M=M_c$, after the time rescaling, equation (\ref{eq:SDE_N}) may be rewritten as
\begin{equation}\label{eq:SDE_N_crit}
	\md R_{t^\prime}\,=\,\frac{1}{2}\frac{\md t^\prime}{R_{t^\prime}}\,+\,\md W_{t^\prime},
	\qquad t^\prime=\frac{2\mu t}{N}.
\end{equation}
This equation describes a critical Bessel process with {\em entrance boundary} at $R=0$, in this case $R_{t^\prime}$ a.s. never hits zero, although it visits its arbitrarily small neighborhood. In this limit the particle system remains stochastic even as $N\to\infty$ and does not approximate the deterministic K-S equations. As $t^\prime\to\infty$, the particle system approaches blow-up infinitely often, but never actually coalesces into a single particle, cf \cite{Bla08b}.

\begin{figure}[h]
\scalebox{0.4}{\includegraphics{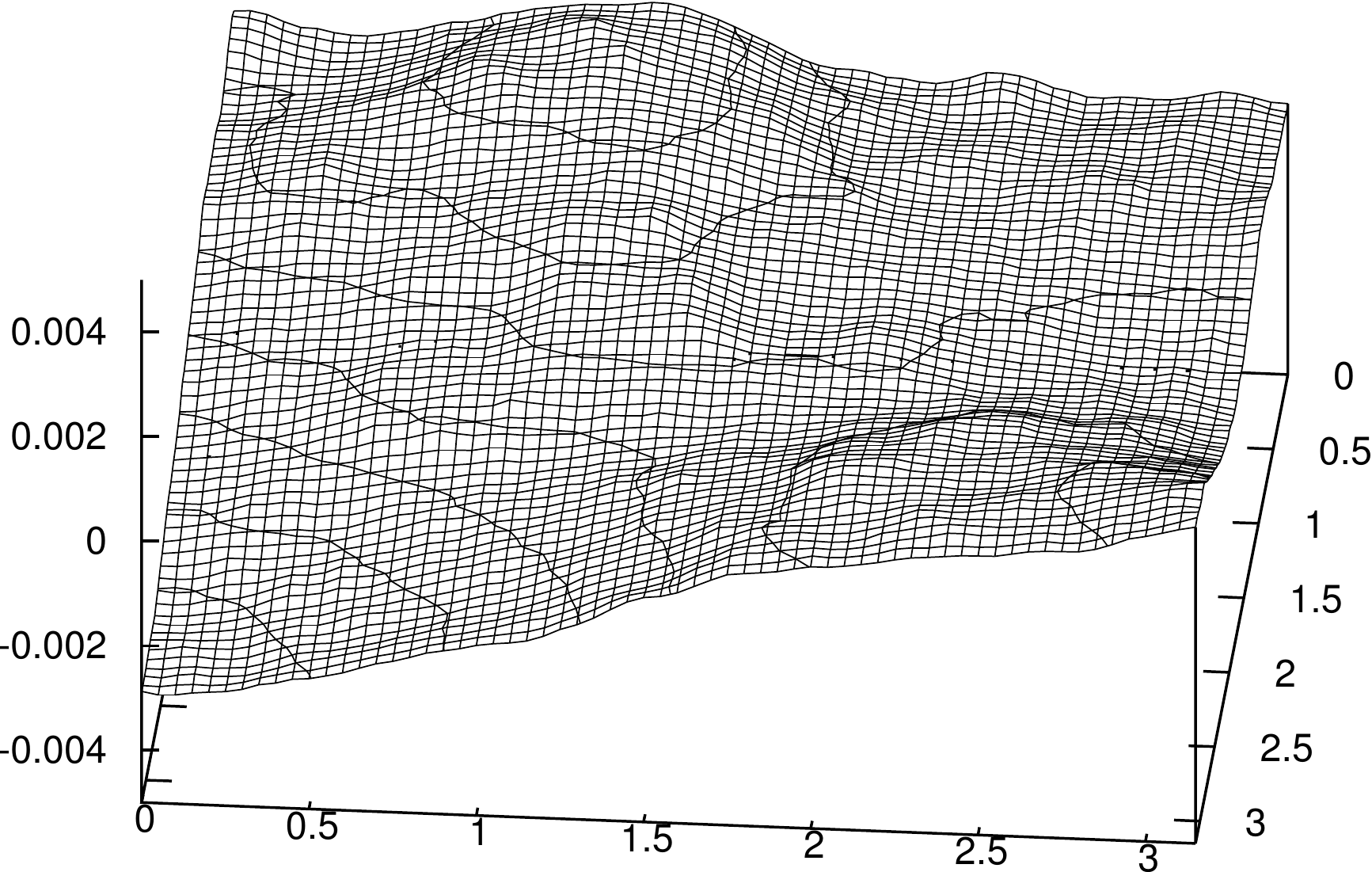}}\hfill\scalebox{0.4}{\includegraphics{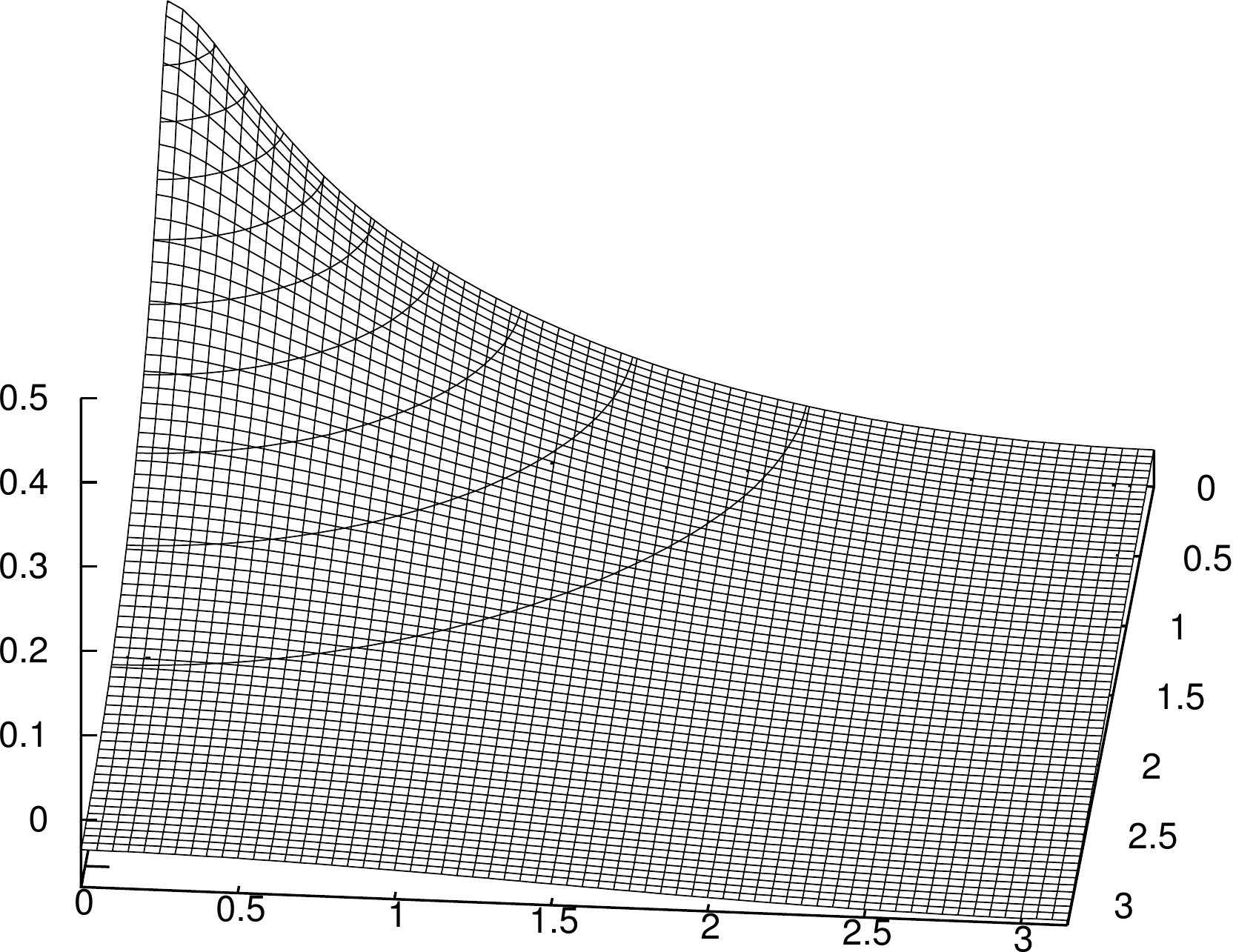}}\\
\caption{\label{fig:supsub}\small Snapshots ($t=10000$) of the concentration field $c(\vv x)-M/(kL)^2$ demonstrating sub- and supercritical behaviors in the elliptic ($\alpha=0$) model. (The constant is subtracted because $c(\vv x)$ itself becomes unbounded in a finite domain as $k\to0$.) The other parameters are as follows: $k=0.01$, $\chi=0.1$, $\mu=0.005$; $\Delta t=0.01$, $\Delta x=0.05$, $4\cdot10^3$ particles. Initially, all particles are distributed randomly with a slight bias towards $(0,0)$. In the supercritical case ($M=0.35$, right plot), a trap is created, while in the subcritical case ($M=0.34$, left plot), all particles eventually spread uniformly over the entire domain. The theoretical critical mass for the elliptic model  with $k=0$ in the entire space is $M_c=2\pi/5\approx1.26$. Notice that only a quarter of the critical mass is needed for a trap in the corner.}
\end{figure}

\paragraph{Formation of singularities.} Our first numerical experiment is designed to test how well our numerical method predicts the value of the critical mass $M_c$ required to develop a singularity from smooth initial data. Notice, first of all, that due to reflection principle for the problem with Neumann boundary conditions, singularity in a corner of the domain, $\Omega=(0,L)^2$, only requires a quarter of the critical mass.  Since the singularity inside of the domain eventually migrates into a corner, we perform our experiment in the corner to begin with. We scatter all particles over $\Omega$ with a slight bias towards $(0,0)$ and vary their total mass observing whether the particles remain aggregated, or spread uniformly over the entire domain. The results are presented in Figure~\ref{fig:supsub}. For the specified values of the parameters, the calculated numerical value of the critical mass lies between 0.34 and 0.35, while the theoretical prediction for elliptic model with $k=0$ in the entire plane is $M_c/4=\pi/10$ --- slightly smaller. The mismatch is a consequence of the simulation in a finite domain, the limit $L\to\infty$, $k\to0$ is non-trivial because equation (\ref{eq:KSc}) with Neumann boundary conditions is not well-posed in $\Omega$ when $\alpha=k=0$. In particular, the boundary of the domain far from the corner where the singularity is formed pulls the particles away from the singularity, effectively increasing the critical mass. 

Interesting phenomena occur when a singularity with mass $M\in(0,M^*_c)$ already exists in the system. In this case the trap absorbs the particles, though the latter may still escape back into the regular component of the particle density. The diffusion process underlying the K-S equations is not uniquely defined by its generator and additional rules for behavior of the particles at singular points must be specified. These rules are not inherently encoded in the K-S equations and are related to non-uniqueness of K-S regularizations.  From the modeling perspective, the exchange of mass between the regular and singular components of the particle density in this regime strongly depends on particular details of the numerical method. Note, however, that if the initial particle density is regular, the smallest possible mass for singularity is greater than $M_c=2M_c^*$,  so this scenario does not occur.
 
\rule{0pt}{4ex}
\begin{figure}[h]
\center\scalebox{0.225}{\includegraphics{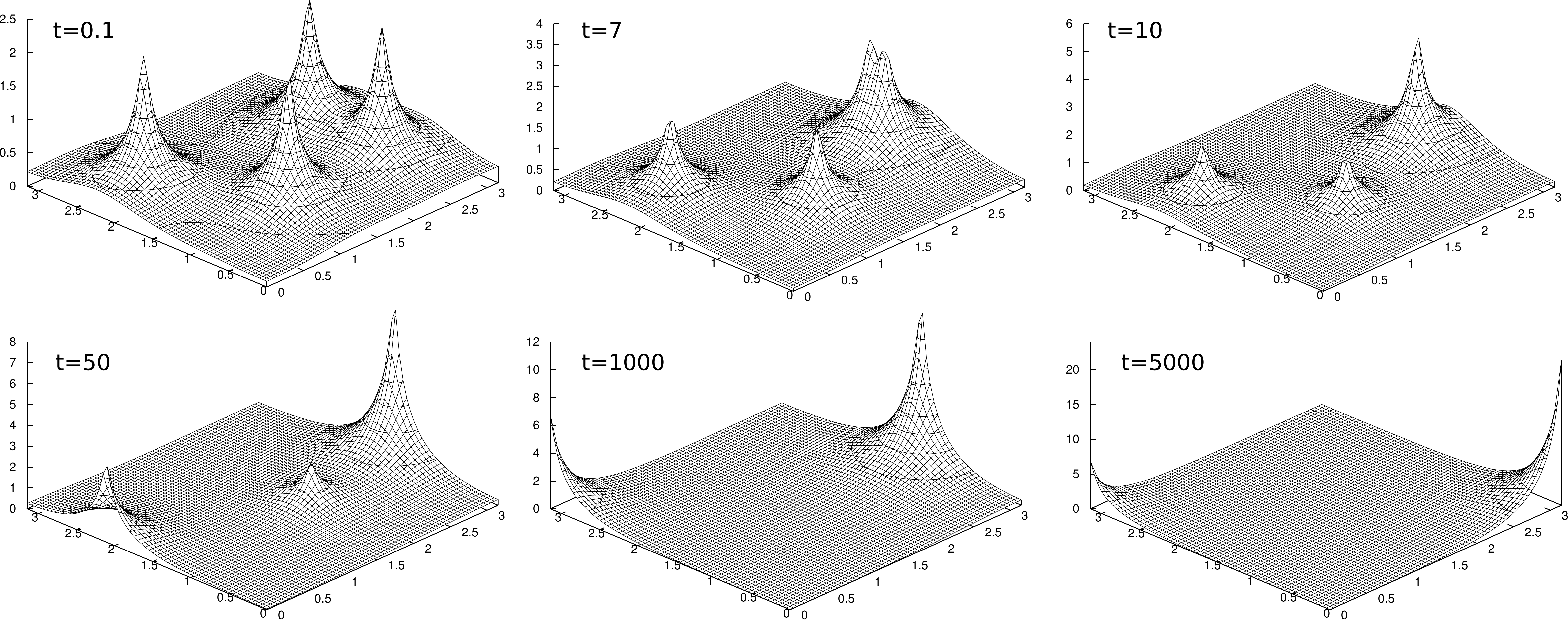}}\\
{\caption{\label{fig:evolv}\small Snapshots of the concentration field $c(\vv x)$ illustrating the motion and interaction of singularities in the elliptic ($\alpha=0$) model. The other parameters are as follows: $k=1$, $\chi=0.1$, $\mu=0.005$, $M=4$; $\Delta t=0.01$, $4\cdot10^3$ particles. Initially, the particles are placed around four distinct locations; as time goes on, the singularities are formed, then they merge and travel towards the boundary of the domain, and finally stabilize in the corners. (Notice the scale change for the $z$-axis.)}\rule{0ex}{1ex}}
\end{figure}

\paragraph{Interaction of singularities} is illustrated in Figure~\ref{fig:evolv}. The singularities attract each other and are also attracted by the boundary of the domain. Equations governing their dynamics may be derived following \cite{Vel04a,Vel04b,Dol09}. If the density $\rho(\vv x)$ is purely atomic, we obtain
\begin{equation}\label{eq:sing_ODE}
  \dot{\vv x}_i(t)\,=\,-\chi\,\frac{\pd}{\pd\!\vv x_i}\sum_{j\neq i}m_j\,V_2(\vv x_i,\vv x_j).
\end{equation}
Here $m_j$-s are the masses of singularities located at $\vv x_j$-s, $V_2(\vv x,\vv y)$ is the Green's function. Note that unlike in the method of Ha\v skovec and Schmeiser \cite{Has10}, in our method, the dynamics of singularities is not imposed by the numerical method explicitly, i.e., it is a natural consequence of the stochastic particle dynamics.

Equations (\ref{eq:sing_ODE}) are identical to equations (\ref{eq:KS_SDE}) of the K-S model with particle diffusivity $\mu$ set to zero. The naive reasoning is that once the traps are formed, the particles' diffusion is dominated by the (infinitely strong) drift and ceases to contribute into dynamics: the stochastic dynamics self-averages and the martingale component does not contribute into the mean drift. In order to verify how well our method approximates equations (\ref{eq:sing_ODE}), we perform a numerical experiment illustrated in Figure~\ref{fig:coll}. In this experiment we create two singularities far from the boundaries in a sufficiently large domain, so that the interaction potential is well-approximated by the fundamental solution (\ref{eq:fs_2d}). We run the K-S solver and track locations of these singularities. Finally, we compare them with locations of the particles evolving according to equation (\ref{eq:sing_ODE}) from the same initial data. A perfect match indicates that the scheme is very successful in dealing with this kind of phenomena. Similarly, in a bounded domain, the particle interaction potential, $V_2(\vv x_i-\vv x_j)$, should be replaced by the Green's function which also encodes interaction with the boundary. In particular, this explains attraction of the singularities to the boundary and the corners of the domain.

\begin{figure}[h]
\rule{0ex}{4ex}\scalebox{0.25}{\includegraphics{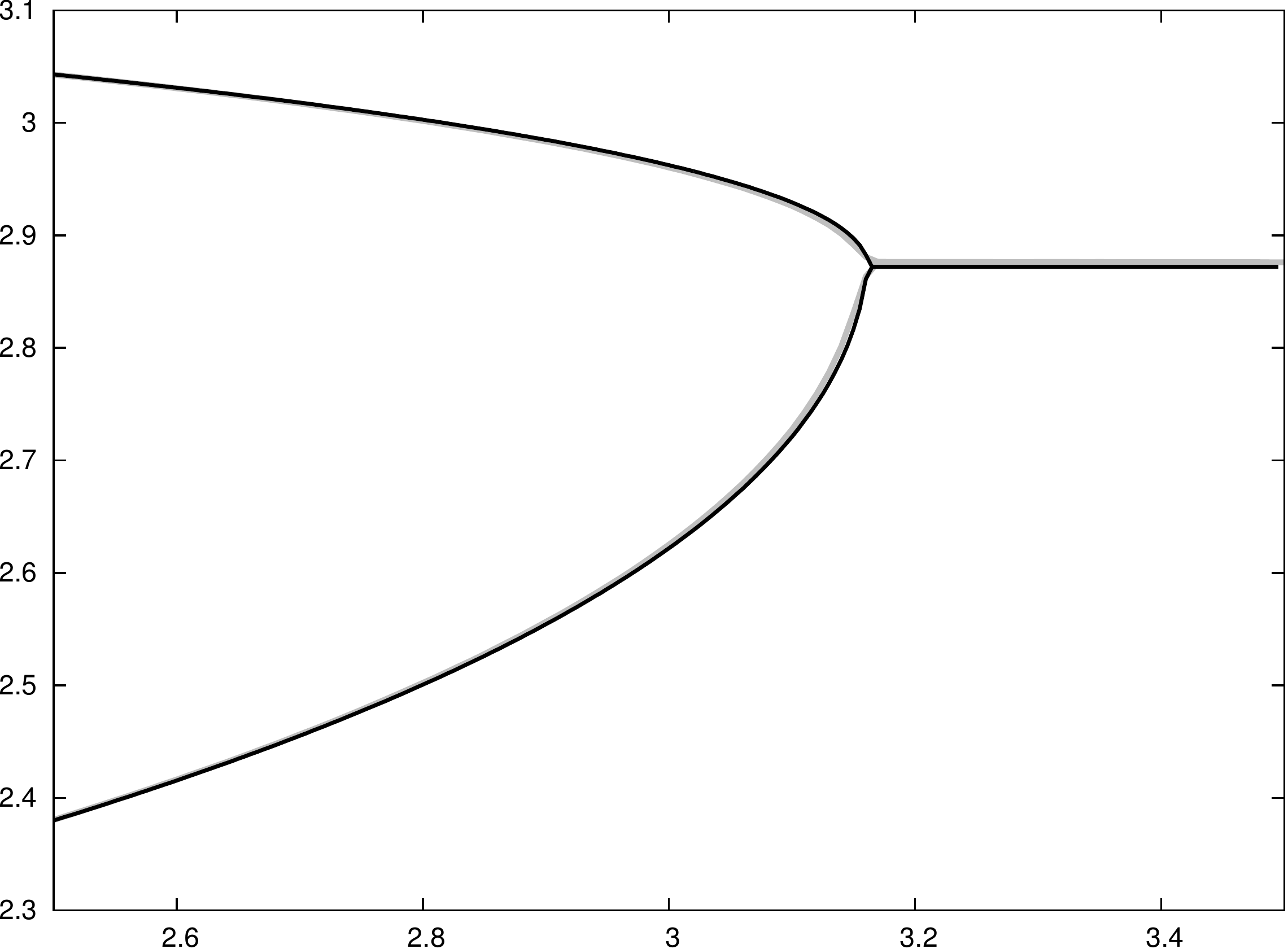}}\hfill
\parbox[b]{3.25in}{\caption{\label{fig:coll}\small Interaction of singularities in the elliptic ($\alpha=0$) model. The particle density $\rho(\vv x)$ is initialized with two delta-functions with masses 6.25 and 18.75. These singularities attract each other and eventually merge; their locations are plotted as functions of time $t$. The thin black line corresponds to the simulation of equation (\ref{eq:sing_ODE}), while the thick grey line --- to simulation of the elliptic Keller-Segel model (the lines are barely distinguishable); $\Delta t=0.01$, $4\cdot10^5$ particles are used. The computational domain is chosen sufficiently large so that the boundary effects are negligible on the presented time scale.}}
\end{figure}


\section{Discussion}

Our investigation of the K-S equations was motivated by the underlying stochastic particle dynamics. A numerical method based on similar ideas has been recently presented by Ha\v skovec and Schmeiser \cite{Has10}.  The greatest advantage of the H-S method is that the singularities are represented as deterministic particles rather than as a tight cluster of ``elemetary'' stochastic particles. This, however, requires an explicit knowledge of evolution equations governing the interaction of singularities (at the very least the knowledge of Green's function in a given domain) and is only applicable to the elliptic K-S model. Moreover, the H-S method requires a direct simulation of an ensemble of interacting particles and thus its computational cost scales as $N^2$ with respect to the total number of particles. By combining the particle dynamics and the PDE dynamics of the concentration field which mediates the particle interaction (the PIC ideas) we are able to significantly reduce this computational cost and deal with a much larger number of interacting particles. We can also treat both parabolic and elliptic models with equal ease. The ultimate numerical method should combine ideas presented in this work, the deterministic particles of the H-S method, and also high order conventional schemes in the regions where the solution remains regular. Such a multi-model implementation remains a subject for the future work.

Even though the non-uniqueness issue is well-known, the details of how the exit-entrance conditions for the underlying diffusion process are related to the exchange of mass between the singularity and the regular part  are not yet well understood. From the modeling perspective, the implication is that a particular numerical scheme provides a regularization of some sort which affects the critical mass $M_c^*$

Another interesting class of stochastic particle models related to our studies of the K-S equations arises in the context of the so-called self-gravitating Brownian particles, see e.g., \cite{Chav07,Chav07b,Sir10}. We suggest a model which fits into this class and bears an intrinsic relation to the K-S model. Consider an ensemble of particles in $d$ dimensions characterized by their masses $m_n$ and evolving according to the stochastic differential equations
\begin{equation}\label{eq:aggr}
  \md\vv X^{(n)}_t\,=\,-\chi\,\frac{\pd}{\pd\!\vv X^{(n)}}\sum_{k\neq n}m_kV_d\big(\vv X^{(n)}_t\!,\,\vv X^{(k)}_t\big)\,\md t\,+\,\sqrt{2\mu M/N_0m_n}\,\md\vv W^{(n)}_t.
\end{equation}
Here $N_0$ is the initial number of particles and $M=\sum m_n$ is their total mass. If initially $m_n=M/N_0$, equation (\ref{eq:aggr}) is precisely (\ref{eq:KS_SDE_n}). The  interaction allows for particle collisions, thus equation (\ref{eq:aggr}) is only valid until the first collision, at which point the colliding particles coalesce into a single particle which acquires their combined mass. Dynamics is then restarted with the remaining particles. In the limit as $N_0\to\infty$, whenever a particle accumulates an $\bO(1)$ mass, its diffusivity becomes zero, i.e., it obeys the deterministic equation (\ref{eq:sing_ODE}) which governs evolution of a point singularity in the K-S model. Therefore we conjecture that in a proper hydrodynamic limit this system is equivalent to the elliptic K-S model. Investigation of this system will be conducted elsewhere.

\bibliographystyle{plain}
\bibliography{../bibliography/bibl}

\end{document}